\documentclass{article}

\usepackage{wrapfig}
\usepackage[preprint]{neurips_2025}

\workshoptitle{Differentiable Systems and Scientific Machine Learning}

\usepackage{amsmath}
\usepackage[utf8]{inputenc} 
\usepackage[T1]{fontenc}    
\usepackage{hyperref}       
\usepackage{url}            
\usepackage{booktabs}       
\usepackage{amsfonts}       
\usepackage{nicefrac}       
\usepackage{microtype}      
\usepackage{xcolor}         
\usepackage{graphicx}
\usepackage{enumitem}
\usepackage{algorithm}   
\usepackage{algpseudocode} 
\usepackage[table]{xcolor}
\definecolor{lightgray}{gray}{0.9}  

\providecommand{\R}{{\mathbb R}}
\DeclareMathOperator*{\argmin}{arg\,min}

\title{Graph Neural Regularizers for PDE Inverse Problems} 

\author{%
  William Lauga \\
  Department of Applied Mathematics\\
  University of Cambridge\\
  wlauga@gmail.com\\
  \And
  James Rowbottom \\
  Department of Applied Mathematics\\
  University of Cambridge\\
  jr908@cam.ac.uk\\
  \And
  Alexander Denker\\
  Department of Computer Science\\
  University College London\\
  a.denker@ucl.ac.uk\\
  \And
  \v{Z}eljko Kereta\\
  Department of Computer Science\\
  University College London\\
  z.kereta@ucl.ac.uk\\
  \And
  Moshe Eliasof\\
  Department of Applied Mathematics\\
  University of Cambridge\\
  me532@cam.ac.uk\\
  \And
  Carola-Bibiane Schönlieb\\
  Department of Applied Mathematics\\
  University of Cambridge\\
  cbs31@cam.ac.uk\\
}

\begin{document}

\maketitle

\begin{abstract}
We present a framework for solving a broad class of ill-posed inverse problems governed by partial differential equations (PDEs), where the target coefficients of the forward operator are recovered through an iterative regularization scheme that alternates between FEM-based inversion and learned graph neural regularization. The forward problem is numerically solved using the finite element method (FEM), enabling applicability to a wide range of geometries and PDEs. By leveraging the graph structure inherent to FEM discretizations, we employ physics-inspired graph neural networks as learned regularizers, providing a robust, interpretable, and generalizable alternative to standard  approaches. Numerical experiments demonstrate that our framework outperforms classical regularization techniques and achieves accurate reconstructions even in highly ill-posed scenarios.
\end{abstract}

\section{Introduction}
Inverse problems governed by partial differential equations (PDEs) arise in many scientific fields, including optical tomography \cite{laiInverseProblemsStationary2018}, seismic imaging \cite{verschuurSeismicDataAnalysis2011}, and combinatorial problems \cite{eliasof2024graph}. In such inverse problems one seeks to recover unknown parameters of a physical system from incomplete or noisy measurements given by the PDE solution. These problems are typically ill-posed, meaning that solutions may not be unique, and can be highly unstable with respect to perturbations in the data. 

In contrast, solving the given forward problem (i.e., the PDE) is typically well-posed and can be accurately discretized using finite element methods (FEM). This produces a mesh-based representation of the PDE domain, which can be interpreted as a graph. This graph structure makes it possible to exploit recent advances in graph neural networks (GNNs) to design learning methods that incorporate the physical structure of the PDE, with applications ranging from medical imaging \cite{zhang2023graph} to geophysics \cite{van2020automated}. 
In this work, we present a framework that leverages physics-inspired GNNs \cite{chamberlain2021grand, eliasof2021pde} to learn a graph neural regularizer \cite{eliasof_learning_2025} (as a data-driven functional) applied to PDE inverse problems.
Our approach directly operates on the FEM mesh, aligning with the numerical discretization of the forward operator. In contrast, black-box methods \cite{ronnebergerUNetConvolutionalNetworks2015, krizhevskyImageNetClassificationDeep2012} learn mappings directly from data; while powerful, they often lack interpretability and are usually limited to pixel-grid scenarios. We discuss related work in Appendix~\ref{app:related_work}. The main contributions of this work are:
\begin{itemize}[itemsep=2pt, topsep=4pt, parsep=0pt, partopsep=0pt, leftmargin=*]
    \item \textbf{Graph-based regularization for PDE inverse problems.}
    We propose a framework that learns \emph{graph neural regularizers} directly on FEM meshes, exploiting their natural graph structure without grid interpolation, and allowing for transfer to irregular domains.
    \item \textbf{Physics-inspired variational formulation.}
    We embed physics-based GNNs as \emph{learned regularizers} within a classical variational inversion scheme, ensuring interpretability and physical consistency.
    \item \textbf{Comprehensive empirical validation.}
    We demonstrate the effectiveness of our framework by applying it to three PDE inverse problems, and show that it can significantly outperform classical methods and is competitive with state-of-the-art black-box methods, even with far fewer parameters.
    
\end{itemize}

\section{Solving PDE Inverse Problems with Graph Neural Regularizers}\label{sec:actual_method} 
We focus on inverse problems governed by a PDE in a domain $D$, for a source term $h$, written as  
\begin{align}
    \label{eq:general_pde}
    \mathcal{L}(u; x) = h \quad \text{in } D,
\end{align}
with appropriate boundary conditions. The operator $\mathcal{L}$ depends on system parameters $x$. The goal is to recover parameters $x$ from noisy and incomplete observations of the solution $u$.

Let $B$ be some observation operator. Then, given observations $y = B(u)$ the inverse problem can be formulated as: recover $x$ such that $B(\mathcal{G}^\dagger(x)) \approx y$, where $\mathcal{G}^\dagger:x \mapsto u$ is the forward operator solving the PDE for $u$ given parameters $x$.
A standard approach for solving this inverse problem is to consider the following variational problem 
\begin{align}
\label{eq:main_ip}
    \hat{x} = \argmin_x \frac{1}{2} \| B(\mathcal{G}^\dagger(x)) - y \|^2 + \alpha \mathcal{R}(x), 
\end{align}
which balances a data-fidelity term with a regularization function $\mathcal{R}(x)$ \cite{noauthor_variational_2009}. A common classical approach that encourages smoothness is Laplacian regularization, given by $\mathcal{R}(x) = \frac{1}{2}x^TLx$, where $L$ is the unnormalized graph Laplacian matrix \cite{nadlerStatisticalAnalysisSemiSupervised2009}. Instead of using a hand-crafted model, we aim to implement the regularizer as a neural network and learn its parameters. One approach to training the regularizer is bilevel learning \cite{haber_learning_2003,kunisch_bilevel_2013}, which adapts the parameters of the network such that the reconstruction in Equation \eqref{eq:main_ip} minimizes some loss. 

We follow the GRIP framework \cite{eliasofDRIPDeepRegularizers2023,eliasof_learning_2025} and consider a lifted version of Equation \eqref{eq:main_ip}, i.e.,
\begin{align}
    \label{eq:lifted_version}
    \argmin_{x_1, \dots, x_L} \frac{1}{2} \| B(\mathcal{G}^\dagger(x_L)) - y \|^2 + 
    \underbrace{
\sum_{l=1}^{L-1} \frac{1}{2} \|x_{l+1} - x_l\|^2
}_{\text{kinetic}}+
\underbrace{\sum_{l=1}^{L-1} \Phi(x_l, f_M, \mathcal{G}, \theta_l)}_{\text{potential}},
\end{align}
with $\hat{x}=x_L$ as the reconstruction and the regularizer given by the kinetic and potential terms above. The potential is given by $\Phi$ and depends on the graph $\mathcal{G}$, meta-data $f_M$ and trainable parameters $\theta_l$ for $l=1,\dots,L-1$. Given a paired dataset $\{ x^{(i)}, y^{(i)} \}_{i=1}^n$ we consider the mean-squared-error loss  
\begin{equation}
\label{eq:loss}
    \min_\theta \sum_{i=1}^n \| \hat{x}(y^{(i)}) - x^{(i)} \|^2,
\end{equation} 
where $\hat{x}(y^{(i)})$ is the reconstruction given by the lifted formulation in Equation \eqref{eq:lifted_version}. 
For optimization of the parameters, we make use of \textit{unrolling} \cite{monga_algorithm_2021}. Here, we use an iterative method to optimize Equation \eqref{eq:lifted_version} for a fixed number of steps and rely on backpropagation to obtain the gradients. In particular, differentiating the regularizer with respect to $x = (x_1, ..., x_L)$ gives the update rule 
\begin{equation}
\label{eq:gnn_step}
x_{l+1} = 2x_l - x_{l-1} + \nabla_{x_l}\Phi(x_l, f_M, \mathcal{G}, \theta_l),
\end{equation}
    where we directly parametrize the gradient $\nabla_{x_l}\Phi$ as a GNN layer \footnote{In practice, we get better results by using the simpler residual update $x_{l+1} = x_l + \nabla_{x_l}\Phi(x_l, f_M, \mathcal{G}, \theta_l).$}. We then solve Equation \eqref{eq:lifted_version} by alternating between regularization steps obtained via a GNN, as in Equation \eqref{eq:gnn_step}, and minimizing the data-fidelity term, using Conjugate Gradient Least-Squares (CGLS). For the meta-data $f_M$ we include a positional encodings, see Appendix \ref{app:meta-data} and Algorithms \ref{alg:training}, \ref{alg:reconstruction}.






\subsection{GNN Architectures} The connection between GNNs and differential equations is well understood \cite{hanContinuousDynamicsGraph2023}. Using physics-inspired GNN architectures ensures alignment between the model and the problem, providing a strong rationale for the approach, as we are applying this framework to physical scenarios. 


\textbf{Graph Neural Diffusion (GRAND).} GRAND \cite{chamberlain2021grand} is an approach that treats GNNs as discretizations of a PDE; it can be seen as a learnable graph-diffusion process. The transformer attention allows it to learn an anisotropic diffusion equation so node features evolve into community patterns, similar to many of the underlying distributions we are trying to recover, thus serving as a useful learnable prior

\textbf{Allen-Cahn Message Passing (ACMP).} The ACMP architecture \cite{wangACMPAllenCahnMessage2022} models the GNN as a particle system with interactive forces and the Allen-Cahn double-well potential. ACMP is an example of a more general class of reaction-diffusion GNNs \cite{choi2023gread,eliasof_graph_2024}. 
ACMP was chosen as a GNN as the double-well encourages node clustering, which is desirable if the target data is known to be, for example, piecewise constant. Our experiments confirm this useful property as ACMP strongly outperforms GRAND in piecewise-constant scenarios, but performance is more similar in the smoother scenarios.
More details on the architectures of GRAND and ACMP can be found in Appendix \ref{app:model_details}.

\section{Experiments and Results}
\label{sec:experiments}
We apply our framework to three representative PDE inverse problems with different geometries to assess its effectiveness and generality. The goal is to evaluate whether our approach can (i) outperform standard GNN baselines, (ii) remain competitive with or even outperform state-of-the-art deep learning methods while using orders of magnitude fewer parameters and less training data, and (iii) provide better interpretability. Moreover, we seek to demonstrate that the framework naturally handles PDEs on irregular meshes without requiring interpolation, making it well-suited to practical inverse problem settings.

We compare against two standard fully-learned baselines: (i) a Graph Convolutional Network (GCN) \cite{kipfSemiSupervisedClassificationGraph2017} and (ii) a U-Net \cite{ronnebergerUNetConvolutionalNetworks2015}. Note, that the U-Net requires an final interpolation step from a structured regular grid to the mesh. 

\textbf{Poisson.}
We consider steady-state diffusion modeled by Poisson's equation $\nabla^2 u = a$ with boundary condition $u|_{\partial D} = 0$. The corresponding inverse problem is to reconstruct the source $a$ from partial measurements of the concentration $u$.
We consider two setups 'dense' and 'sparse', depending on the number of observations of $u$. Results are summarized in Table \ref{tab:poisson}.

\begin{table}[t]
  \centering
  \caption{Results for the Poisson PDE ($\pm$ absolute standard deviation). We considered a dense sampling setup with 60\% of vertices observed, and a sparse setup with 10\% of vertices observed. We report the
  MSE between reconstructed and ground-truth coefficients, and the data-fit measuring agreement with observed data.}
  \label{tab:poisson}
  \begin{tabular}{lcccc}
    \toprule
    & \multicolumn{2}{c}{Poisson (Dense)} & \multicolumn{2}{c}{Poisson (Sparse)} \\
    \cmidrule(lr){2-3} \cmidrule(lr){4-5}
    Method & MSE ($\downarrow$) & Data-fit ($\downarrow$) & MSE ($\downarrow$) & Data-fit ($\downarrow$) \\
    \midrule
    Laplacian & $0.1100~(\pm 0.035)$ & \cellcolor{lightgray}$0.0023~(\pm 0.001)$ & $0.2201~(\pm 0.074)$ & \cellcolor{lightgray}$0.0067~(\pm 0.007)$ \\
    U-Net      & $0.1220~(\pm 0.110)$ & $0.0111~(\pm 0.013)$ & $0.1587~(\pm 0.100)$ & $0.0905~(\pm 0.032)$ \\
    GCN       & $0.0646~(\pm 0.031)$ & $0.0155~(\pm 0.024)$ & $0.1879~(\pm 0.130)$ & $0.1207~(\pm 0.110)$ \\
    \midrule
    \textbf{Our framework} \\
    GRIP (GRAND) & $0.0647~(\pm 0.025)$ & $0.0259~(\pm 0.050)$ & \cellcolor{lightgray}$0.1441~(\pm 0.066)$ & $0.0835~(\pm 0.020)$ \\
    GRIP (ACMP) & \cellcolor{lightgray}$0.0586~(\pm 0.023)$ & $0.0108~(\pm 0.013)$ & $0.1443~(\pm 0.063)$ & $0.0826~(\pm 0.036)$ \\
    \bottomrule
  \end{tabular}
\end{table}

\begin{figure}[t]
  \centering
\includegraphics[width=1.0\textwidth]{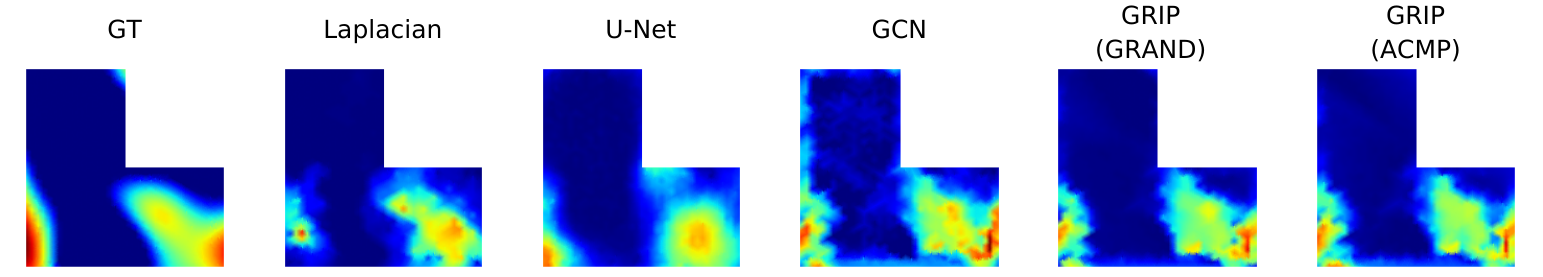}
  \vspace{-0.5cm}
  \caption{Poisson (sparse) reconstruction examples}
\end{figure}

\textbf{Inverse Wave Scattering.}
We seek to infer the material properties of a medium by wave propagation in the frequency domain as described in \cite{molinaro_neural_nodate}. The governing PDE is the Helmholtz equation
\begin{equation}
    -\Delta u - \omega^{2} au = 0, \quad \text{in }  D, \text{ and } u = g, \quad \text{on } \partial D,
\end{equation}
where $\omega$ is a frequency and we have the Dirichlet boundary condition $g\in H^\frac{1}{2}(\partial D)$. We observe data on the boundary $\partial D$ with a Dirichlet-to-Neumann map, and we try to solve the linearized version of the problem. Details of the linearization and of the boundary conditions are given in Appendix~\ref{app:pde_details}. We study this problem, whose results are summarized in Table \ref{tab:inverse_wave_eit}, on the square domain \(D = [0,1]^2\). The coefficient $a$ is sampled from a distribution containing 1-4 square-shaped inclusions that are randomly distributed across the domain. 

\begin{figure}[t]
  \centering
  \includegraphics[width=0.95\textwidth]{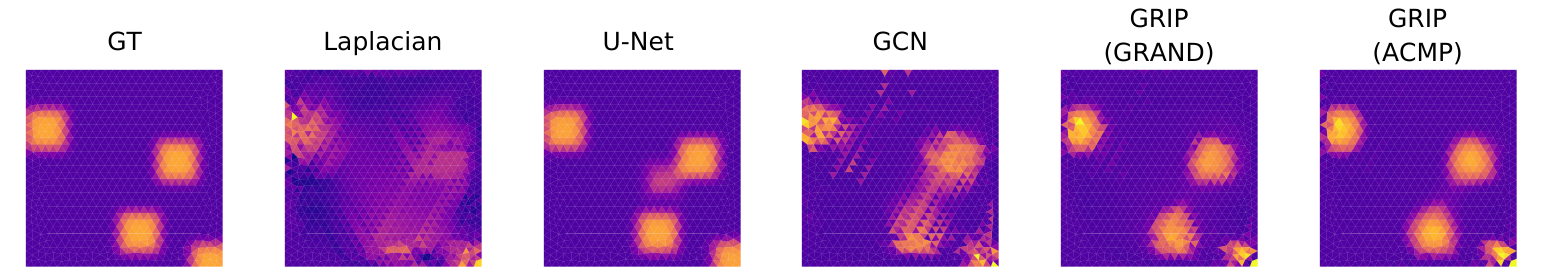}
  \vspace{-0.2cm}
  \caption{Inverse Wave Scattering reconstruction examples}
\end{figure}

\begin{figure}[t]
  \centering
\includegraphics[width=0.95\textwidth]{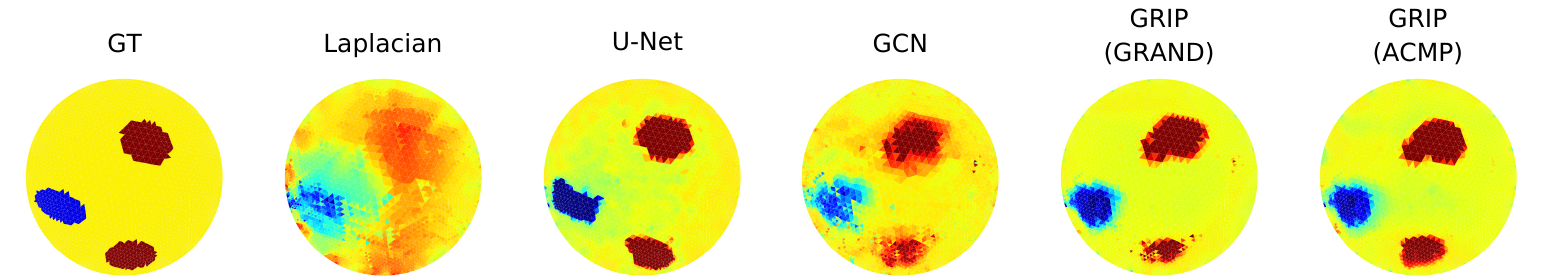}
  \vspace{-0.2cm}
  \caption{EIT reconstruction examples}
\end{figure}

\textbf{Electrical Impedance Tomography.}
EIT is a noninvasive imaging technique that uses electrical measurements taken on the boundary of an object to reconstruct the conductivity of its interior \cite{borcea_electrical_2002}.
We have a domain $D$, and $U^{L}_{l=1} e_l$ is the set of $L$ electrodes attached to the boundary $\partial D$. The electric potential $u$ is given by
\begin{equation}
    - \nabla \cdot (\sigma \nabla u) = 0 \quad \text{in}\quad D ,
\end{equation}
where $\sigma$ is the conductivity distribution. Realistic boundary conditions are given by the Complete Electrode Model (CEM), as described in \cite{somersalo_existence_1992}. The goal is to reconstruct $\sigma$ given electrode measurements. We linearize the forward operator around a constant background conductivity. Results are reported in Table~\ref{tab:inverse_wave_eit}. Here, we also compare to a state-of-the-art U-Net from \cite{denker_deep_2025}, which has about $10^7$ parameters, whereas the GNNs in our framework have < 1,500 parameters.
\vspace{-0.2cm}

\begin{table}[hb]
\centering
\caption{Results for Inverse Wave Scattering and EIT ($\pm$ absolute standard deviation). We report the MSE between reconstructed and ground-truth coefficients/conductivity, and the data-fit measuring agreement with observed boundary or electrode data.}
\label{tab:inverse_wave_eit}
\begin{tabular}{lcccc}
\toprule
& \multicolumn{2}{c}{Inverse Wave Scattering} & \multicolumn{2}{c}{EIT} \\
\cmidrule(lr){2-3} \cmidrule(lr){4-5}
Method & MSE ($\downarrow$) & Data-fit ($\downarrow$) & MSE ($\downarrow$) & Data-fit ($\downarrow$) \\
\midrule
Laplacian & $0.3700~(\pm 0.083)$ & $0.0047~(\pm 0.027)$ & $0.6967~(\pm 0.160)$ & $0.0675~(\pm 0.091)$ \\
U-Net & \cellcolor{lightgray}$0.0730~(\pm 0.091)$ & $0.0506~(\pm 0.023)$ & \cellcolor{lightgray}$0.3009~(\pm 0.120)$ & \cellcolor{lightgray}$0.0974~(\pm 0.076)$ \\
GCN & $0.2402~(\pm 0.065)$ & $0.0601~(\pm 0.023)$ & $0.3510~(\pm 0.065)$ & $0.3584~(\pm 0.100)$ \\
\midrule
\textbf{Our framework} \\
GRIP (GRAND) & $0.0953~(\pm 0.036)$ & $0.0022~(\pm 0.001)$ & $0.3887~(\pm 0.16)$ & $0.2229~(\pm 0.51)$ \\
GRIP (ACMP) & $0.0769~(\pm 0.047)$ & \cellcolor{lightgray}$0.0020~(\pm 0.001)$ & $0.3471~(\pm 0.16)$ & $0.2199~(\pm 0.53)$ \\
\bottomrule
\end{tabular}
\end{table}

\section{Conclusion}
We introduced a general framework for solving PDE-based inverse problems by learning graph-based regularizers with physics-inspired GNN architectures that exploit the natural graph structure of FEM discretizations. Across three representative PDE inverse problems we demonstrated that our framework consistently outperforms classical variational regularization methods and achieves performance competitive with state-of-the-art black-box deep learning approaches, occasionally even outperforming them as in the Poisson PDE example.

An important advantage of our method is that it preserves interpretability by incorporating the regularizer into a variational framework. The forward operator and data-fidelity term remain explicit, unlike purely supervised mappings such as U-Net. In contrast to most convolutional neural networks, such as the U-Net, our GNN-based regularizer operates natively on meshes and can be flexibly combined with iterative solvers such as CGLS, without intermediate interpolation steps.
Further work will include combining our approach with Gauss-Newton iterations, to mitigate linearization errors, and extending it to more complex geometries.




\bibliography{references, bib, bib2}
\bibliographystyle{abbrv}

\clearpage

\appendix

\section{Related Work}
\label{app:related_work}
\textbf{On Learned Regularizers.}
Recent studies have explored learning the regularizer $\mathcal{R}$ from data, see e.g. \cite{habring_neuralnetworkbased_2024,mukherjee_learned_2023, mardaniRecurrentGenerativeAdversarial2017,hertrich_learning_2025}. These data-driven regularizers often demonstrate improved performance compared to traditional handcrafted penalties. More recently, \cite{eliasof_learning_2025} introduced a graph-based learned regularization framework, demonstrating that incorporating structural relationships enables the method to generalize effectively to more complex and irregular data topologies.

\textbf{On Neural Operators.} 
Neural operators aim to learn mappings between function spaces \cite{nikola_kovachki_neural_2023}. There are two main design philosophies: (i) transforming functions into a different basis, such as in the Fourier Neural Operator \cite{zongyi_li_fourier_2021}, and learning the mapping there, or (ii) operating directly on irregular domains through graph-based representations. We make use of the latter approach and implement the mapping using GNNs. While GNNs have been used in EIT in the literature \cite{herzberg_domain_2023,herzberg_graph_2021}, these works rely on supervised learning, directly learning a mapping from measurements to reconstructions. In contrast, we learn a regularizer to be used in the variational framework.

\section{PDE Descriptions and Details}
\label{app:pde_details}

This section provides additional details on the PDE inverse problems considered in the experiments. We solve the PDEs using FEM, see \cite{brenner_mathematical_2008} for a mathematical overview. In particular, we employ a triangulation with $N$ nodes of the domain $D$. This triangulation defines a graph structure, where each node defines a vertex and we add an edge whenever two nodes are connected. Further, this triangulation defines a finite-dimensional subspace as the set of piecewise linear functions restricted to the triangulation. 

\paragraph{Linearization.}
To simplify training, we solve the linearized inverse problems by assuming a known reference coefficient $a_\mathrm{ref}$ and attempt to reconstruct a perturbation $\partial a$, yielding
\begin{equation}
    \Tilde{F}(a_\mathrm{ref} + \partial a; a_\mathrm{ref}) \approx F(a_\mathrm{ref}) + J_{a_\mathrm{ref}}(\partial a),
\end{equation}
where $J_{a_\mathrm{ref}} = \nabla F(a_\mathrm{ref})$ is the Jacobian of the forward map.

\subsection{Poisson Equation}
\paragraph{Data generation.}
Ground-truth source fields $a(x)$ were randomly sampled at mesh vertices and smoothed using a Gaussian kernel. Forward solutions were obtained via finite element discretization on an L-shaped domain. The L-shaped mesh used had 1990 elements.
Two setups were used:
\begin{itemize}[itemsep=2pt, topsep=2pt, leftmargin=*]
    \item \textbf{Dense:} $60\%$ of vertices observed.
    \item \textbf{Sparse:} $10\%$ of vertices observed.
\end{itemize}

\paragraph{FEM and weak formulation.} For the Poisson PDE the weak formulation reads
\begin{align*}
    \text{Find } u \in H_0^1(D) \text{ such that } \int_D \nabla u \cdot \nabla v dx = - \int_D a v dx, \quad \forall v \in H_0^1(D).
\end{align*}
Let $\{ \phi_j\}_{j=1}^N$ be the FEM basis of piecewise linear hat functions. If we approximate 
\begin{align*}
    u(x) \approx \sum_{j=1}^N u_j \phi_j(x), \quad a(x) \approx \sum_{j=1}^N a_j \phi_j(x),
\end{align*}
we obtain a linear system for the coefficients $\mathbf{u}=(u_1, \dots, u_N)$ as 
\begin{align*}
    K \mathbf{u} = - M \mathbf{a},
\end{align*}
with $\mathbf{a}=(a_1, \dots, a_N)$, stiffness matrix $K \in \R^{N \times N}$ and mass matrix $M \in \R^{N \times N}$ defined by 
\begin{align*}
    K_{i,j}= \int_D \nabla \phi_j \cdot \nabla\phi_j dx, \quad M_{i,j} = \int_D \phi_i \phi_j dx.
\end{align*}
Note that the linear system in linear in $\mathbf{a}$, so we are able to directly use CGLS to approximate the solution.

\subsection{Inverse Wave Scattering (Helmholtz Equation)}

\paragraph{Boundary conditions}

We observe data on the boundary $\partial D$ with the Dirichlet-to-Neumann map,
\begin{equation}
    \Lambda_a : H^{1/2}(\partial D) \;\longrightarrow\;H^{-1/2}(\partial D),
\end{equation}
\begin{equation}
    \Lambda_a[g] = \left.\frac{\partial u}{\partial \nu}\right|_{\partial D}, \quad \forall g \in H^{1/2}(\partial D).
\end{equation}
Our goal is to reconstruct the wave coefficient $a$ from the map $\Lambda_a$.

\paragraph{Coefficient distribution.}
The coefficients are defined as mixtures of fourth-order Gaussian-like inclusions
\begin{equation}
    a(x,y) = \sum_{k=1}^m \exp\!\big(-c(x - c_{1,k})^4 - c(y - c_{2,k})^4\big),
\end{equation}
with $m \in [1,4]$ and $c = 2\times10^{4}/3$.

The distribution models a homogeneous medium featuring 1–4 square inclusions randomly positioned within the domain. For each sampled coefficient field, 20 Dirichlet boundary conditions are imposed with trigonometric boundary data. The mesh used here had 1578 elements.

\paragraph{FEM and weak formulation.} The weak formulation for the Helmholtz equation reads 
\begin{align*}
    \text{Find } u \in H^1(D) \text{ such that } \int_D \nabla u \cdot \nabla v dx - \omega^2 \int_D a u v dx = 0, \quad \forall v \in H_0^1(D),
\end{align*}
with $u_{\vert \partial D}=g$. Further, we assume that both $\omega$ and $a$ are real valued. 

\subsection{Electrical Impedance Tomography}


\paragraph{Dataset generation}
We used synthetically-generated data on a circular domain, where each datapoint has 1-3 ellipses with a random constant conductivity. The mesh used for the EIT discretization had 4984 elements.

\paragraph{FEM and weak formulation.} We employ the FEM method described in \cite{denker_deep_2025}, which augments the standard FEM system to enforce conservation of charge and grounding condition, i.e., the sum of the voltages measured at the electrodes $e_l$ has to sum to zero. We compute the Jacobian used in the linearization by the method described in \cite{polydorides_matlab_2002} to reduce the computational cost.

\section{Model Architecture, Hyperparameters, and Training}
\label{app:model_details}

\subsection{On Meta-Data for PDE Inverse Problems.} In 
\label{app:meta-data}\cite{eliasof_learning_2025} it is shown that including meta-data, denoted $f_M$, can drastically improve model performance. Hence, we concatenate the graph positional-encodings and then apply a 2-layer multi-layer perceptron (MLP) before the regularization steps.

\subsection{GNN Architectures}

\paragraph{GRAND (Graph Neural Diffusion).}
GRAND treats message passing as a diffusion process. If $X$ are some embedded node features, GRAND tries to solve:
\begin{equation}
    X(T) = X(0) + \int_{0}^{T} \frac{\partial X(t)}{\partial t}\, dt,
\end{equation}
where
\begin{equation}
    \frac{\partial X_i(t)}{\partial t} 
= \sum_{j : (i,j) \in E_0} 
a\bigl(X_i(t), X_j(t)\bigr)\,\bigl(X_j(t) - X_i(t)\bigr),
\end{equation}
where $a(\cdot,\cdot)$ is a learnable attention function.
We use the linear variant (GRAND-l) for computational efficiency.
Implicit Euler time discretization is applied.

\paragraph{ACMP (Allen–Cahn Message Passing).}
ACMP models message passing as a reaction–diffusion system:
\begin{equation}
    \frac{\partial}{\partial t} x_i(t) 
= \alpha \odot \sum_{j \in N_i} \bigl(a(x_i(t), x_j(t)) - \beta \bigr)\,\bigl(x_j(t) - x_i(t)\bigr) + \delta \odot x_i(t) \odot \bigl(1 - x_i(t) \odot x_i(t)\bigr),
\end{equation}
where $\beta$ models attractive/repulsive forces between particles, and $x_i(t) \odot \bigl(1 - x_i(t) \odot x_i(t)\bigr)$ is the Allen-Cahn double-well potential. Setting $\beta$ and $\delta$ to 0 reduces ACMP to GRAND. In our experiments we set $\beta$ to 0 as we are only after the power of the double-well, and $\alpha$ and $\delta$ are learnable parameters.

For both GRAND and ACMP we use 32 layers.

\subsection{Training Procedure}

\begin{algorithm}[t]
\caption{Training the Graph Neural Regularizer}
\label{alg:training}
\begin{algorithmic}[1]
\Require Forward operator $\mathbf{G^\dagger}$, dataset $\{(x,y)\}$, number of iterations $N_\text{max}$
\While{not converged}
    \For{($x$,$y$) in dataset}
        \State $z$ = 0
        \For{$k=1$ to $N_{max}$}
            \State Run CGLS steps on $z$
            \State Concatenate $z$ with positional encodings 
            \State Apply 2-layer MLP
            \State Run $\text{GNN}_{\theta}$ on $z$
        \EndFor
        \State unembed $z$
        \State Compute loss of $z$ vs. $x$ as in Equation \eqref{eq:loss}
        \State Take gradient step w.r.t.~$(\theta)$
    \EndFor
\EndWhile
\State \textbf{return} Fitted parameters $(\theta)$
\end{algorithmic}
\end{algorithm}

\begin{algorithm}[t]
\caption{Reconstruction}
\label{alg:reconstruction}
\begin{algorithmic}[1]
\Require Forward operator $\mathbf{G^\dagger}$, observation $y$, number of iterations $N_\text{max}$
\State $z = 0$
\For{$k=1$ to $N_\text{max}$}
    \State Run CGLS steps on $z$
    \State Concatenate $z$ with positional encodings 
    \State Apply 2-layer MLP
    \State Run $\text{GNN}_{\theta}$ on $z$
\EndFor
\State unembed $z$
\State \textbf{return} $z$ (the reconstruction)
\end{algorithmic}
\end{algorithm}

The training and reconstruction algorithms are given by Algorithms \ref{alg:training} and \ref{alg:reconstruction}. We use $15$ unrolling steps, with $15-30$ CGLS iterations per step. Models were trained for $100$ to $400$ epochs with early stopping. For training, we used $100$ to $200$ data points for the GNNs and $500$ for the U-Net. Adam optimizer was used with a learning rate of $10^{-3}$. We evaluated the model on a held-out set of $100$ test data points. The U-Nets had approximately $10^7$ parameters, whereas the GNNs had fewer than 1500.

\end{document}